\documentstyle[amssymb,12pt]{amsart}

\newcommand{\HypFsimple}[2]{\sideset{_#1}{_#2}{F}}
\newcommand{\HypF}[5]{\sideset{_#1}{_#2}{F}\!\left[\matrix  
#3\\#4\endmatrix;{\displaystyle #5}\right]}
\def\po#1#2{(#1)_{#2}}
\newcommand{\floor}[1]{\left\lfloor#1\right\rfloor}

\newtheorem{thm}{Theorem}

\newtheorem{cor}[thm]{Corollary}
\newtheorem{lem}[thm]{Lemma}

\numberwithin{equation}{section}

\catcode`\@=11
\font\tenln    = line10
\font\tenlnw   = linew10

\newskip\Einheit \Einheit=0.5cm
\newcount\xcoord \newcount\ycoord
\newdimen\xdim \newdimen\ydim \newdimen\PfadD@cke \newdimen\Pfadd@cke

\newcount\@tempcnta
\newcount\@tempcntb

\newdimen\@tempdima
\newdimen\@tempdimb

\newdimen\@wholewidth
\newdimen\@halfwidth

\newcount\@xarg
\newcount\@yarg
\newcount\@yyarg
\newbox\@linechar
\newbox\@tempboxa
\newdimen\@linelen
\newdimen\@clnwd
\newdimen\@clnht

\newif\if@negarg

\def\@whilenoop#1{}
\def\@whiledim#1\do #2{\ifdim #1\relax#2\@iwhiledim{#1\relax#2}\fi}
\def\@iwhiledim#1{\ifdim #1\let\@nextwhile=\@iwhiledim
        \else\let\@nextwhile=\@whilenoop\fi\@nextwhile{#1}}

\def\@whileswnoop#1\fi{}
\def\@whilesw#1\fi#2{#1#2\@iwhilesw{#1#2}\fi\fi}
\def\@iwhilesw#1\fi{#1\let\@nextwhile=\@iwhilesw
         \else\let\@nextwhile=\@whileswnoop\fi\@nextwhile{#1}\fi}

\def\thinlines{\let\@linefnt\tenln \let\@circlefnt\tencirc
  \@wholewidth\fontdimen8\tenln \@halfwidth .5\@wholewidth}
\def\thicklines{\let\@linefnt\tenlnw \let\@circlefnt\tencircw
  \@wholewidth\fontdimen8\tenlnw \@halfwidth .5\@wholewidth}
\thinlines

\def\PfadDicke#1{\PfadD@cke#1 \divide\PfadD@cke by2 \Pfadd@cke\PfadD@cke  
\multiply\PfadD@cke by2}
\long\def\LOOP#1\REPEAT{\def\BODY{#1}\ITERATE}
\def\ITERATE{\BODY \let\next\ITERATE \else\let\next\relax\fi \next}
\let\REPEAT=\fi
\def\Punkt{\hbox{\raise-2pt\hbox to0pt{\hss$\ssize\bullet$\hss}}}
\def\DuennPunkt(#1,#2){\unskip
  \raise#2 \Einheit\hbox to0pt{\hskip#1 \Einheit
          \raise-2.5pt\hbox to0pt{\hss$\bullet$\hss}\hss}}
\def\NormalPunkt(#1,#2){\unskip
  \raise#2 \Einheit\hbox to0pt{\hskip#1 \Einheit
          \raise-3pt\hbox to0pt{\hss\twelvepoint$\bullet$\hss}\hss}}
\def\DickPunkt(#1,#2){\unskip
  \raise#2 \Einheit\hbox to0pt{\hskip#1 \Einheit
          \raise-4pt\hbox to0pt{\hss\fourteenpoint$\bullet$\hss}\hss}}
\def\Kreis(#1,#2){\unskip
  \raise#2 \Einheit\hbox to0pt{\hskip#1 \Einheit
          \raise-4pt\hbox to0pt{\hss\fourteenpoint$\circ$\hss}\hss}}

\def\Line@(#1,#2)#3{\@xarg #1\relax \@yarg #2\relax
\@linelen=#3\Einheit
\ifnum\@xarg =0 \@vline
  \else \ifnum\@yarg =0 \@hline \else \@sline\fi
\fi}

\def\@sline{\ifnum\@xarg< 0 \@negargtrue \@xarg -\@xarg \@yyarg -\@yarg
  \else \@negargfalse \@yyarg \@yarg \fi
\ifnum \@yyarg >0 \@tempcnta\@yyarg \else \@tempcnta -\@yyarg \fi
\ifnum\@tempcnta>6 \@badlinearg\@tempcnta0 \fi
\ifnum\@xarg>6 \@badlinearg\@xarg 1 \fi
\setbox\@linechar\hbox{\@linefnt\@getlinechar(\@xarg,\@yyarg)}%
\ifnum \@yarg >0 \let\@upordown\raise \@clnht\z@
   \else\let\@upordown\lower \@clnht \ht\@linechar\fi
\@clnwd=\wd\@linechar
\if@negarg \hskip -\wd\@linechar \def\@tempa{\hskip -2\wd\@linechar}\else
     \let\@tempa\relax \fi
\@whiledim \@clnwd <\@linelen \do
  {\@upordown\@clnht\copy\@linechar
   \@tempa
   \advance\@clnht \ht\@linechar
   \advance\@clnwd \wd\@linechar}%
\advance\@clnht -\ht\@linechar
\advance\@clnwd -\wd\@linechar
\@tempdima\@linelen\advance\@tempdima -\@clnwd
\@tempdimb\@tempdima\advance\@tempdimb -\wd\@linechar
\if@negarg \hskip -\@tempdimb \else \hskip \@tempdimb \fi
\multiply\@tempdima \@m
\@tempcnta \@tempdima \@tempdima \wd\@linechar \divide\@tempcnta \@tempdima
\@tempdima \ht\@linechar \multiply\@tempdima \@tempcnta
\divide\@tempdima \@m
\advance\@clnht \@tempdima
\ifdim \@linelen <\wd\@linechar
   \hskip \wd\@linechar
  \else\@upordown\@clnht\copy\@linechar\fi}

\def\@hline{\ifnum \@xarg <0 \hskip -\@linelen \fi
\vrule height\Pfadd@cke width \@linelen depth\Pfadd@cke
\ifnum \@xarg <0 \hskip -\@linelen \fi}

\def\@getlinechar(#1,#2){\@tempcnta#1\relax\multiply\@tempcnta 8
\advance\@tempcnta -9 \ifnum #2>0 \advance\@tempcnta #2\relax\else
\advance\@tempcnta -#2\relax\advance\@tempcnta 64 \fi
\char\@tempcnta}

\def\Vektor(#1,#2)#3(#4,#5){\unskip\leavevmode
  \xcoord#4\relax \ycoord#5\relax
      \raise\ycoord \Einheit\hbox to0pt{\hskip\xcoord \Einheit
         \Vector@(#1,#2){#3}\hss}}

\def\Vector@(#1,#2)#3{\@xarg #1\relax \@yarg #2\relax
\@tempcnta \ifnum\@xarg<0 -\@xarg\else\@xarg\fi
\ifnum\@tempcnta<5\relax
\@linelen=#3\Einheit
\ifnum\@xarg =0 \@vvector
  \else \ifnum\@yarg =0 \@hvector \else \@svector\fi
\fi
\else\@badlinearg\fi}

\def\@hvector{\@hline\hbox to 0pt{\@linefnt
\ifnum \@xarg <0 \@getlarrow(1,0)\hss\else
    \hss\@getrarrow(1,0)\fi}}

\def\@vvector{\ifnum \@yarg <0 \@downvector \else \@upvector \fi}

\def\@svector{\@sline
\@tempcnta\@yarg \ifnum\@tempcnta <0 \@tempcnta=-\@tempcnta\fi
\ifnum\@tempcnta <5
  \hskip -\wd\@linechar
  \@upordown\@clnht \hbox{\@linefnt  \if@negarg
  \@getlarrow(\@xarg,\@yyarg) \else \@getrarrow(\@xarg,\@yyarg) \fi}%
\else\@badlinearg\fi}

\def\@upline{\hbox to \z@{\hskip -.5\Pfadd@cke \vrule width \Pfadd@cke
   height \@linelen depth \z@\hss}}

\def\@downline{\hbox to \z@{\hskip -.5\Pfadd@cke \vrule width \Pfadd@cke
   height \z@ depth \@linelen \hss}}

\def\@upvector{\@upline\setbox\@tempboxa\hbox{\@linefnt\char'66}\raise
     \@linelen \hbox to\z@{\lower \ht\@tempboxa\box\@tempboxa\hss}}

\def\@downvector{\@downline\lower \@linelen
      \hbox to \z@{\@linefnt\char'77\hss}}

\def\@getlarrow(#1,#2){\ifnum #2 =\z@ \@tempcnta='33\else
\@tempcnta=#1\relax\multiply\@tempcnta \sixt@@n \advance\@tempcnta
-9 \@tempcntb=#2\relax\multiply\@tempcntb \tw@
\ifnum \@tempcntb >0 \advance\@tempcnta \@tempcntb\relax
\else\advance\@tempcnta -\@tempcntb\advance\@tempcnta 64
\fi\fi\char\@tempcnta}

\def\@getrarrow(#1,#2){\@tempcntb=#2\relax
\ifnum\@tempcntb < 0 \@tempcntb=-\@tempcntb\relax\fi
\ifcase \@tempcntb\relax \@tempcnta='55 \or
\ifnum #1<3 \@tempcnta=#1\relax\multiply\@tempcnta
24 \advance\@tempcnta -6 \else \ifnum #1=3 \@tempcnta=49
\else\@tempcnta=58 \fi\fi\or
\ifnum #1<3 \@tempcnta=#1\relax\multiply\@tempcnta
24 \advance\@tempcnta -3 \else \@tempcnta=51\fi\or
\@tempcnta=#1\relax\multiply\@tempcnta
\sixt@@n \advance\@tempcnta -\tw@ \else
\@tempcnta=#1\relax\multiply\@tempcnta
\sixt@@n \advance\@tempcnta 7 \fi\ifnum #2<0 \advance\@tempcnta 64 \fi
\char\@tempcnta}

\def\Diagonale(#1,#2)#3{\unskip\leavevmode
  \xcoord#1\relax \ycoord#2\relax
      \raise\ycoord \Einheit\hbox to0pt{\hskip\xcoord \Einheit
         \Line@(1,1){#3}\hss}}
\def\AntiDiagonale(#1,#2)#3{\unskip\leavevmode
  \xcoord#1\relax \ycoord#2\relax 
      \raise\ycoord \Einheit\hbox to0pt{\hskip\xcoord \Einheit
         \Line@(1,-1){#3}\hss}}
\def\Pfad(#1,#2),#3\endPfad{\unskip\leavevmode
  \xcoord#1 \ycoord#2 \thicklines\ZeichnePfad#3\endPfad\thinlines}
\def\ZeichnePfad#1{\ifx#1\endPfad\let\next\relax
  \else\let\next\ZeichnePfad
    \ifnum#1=1
      \raise\ycoord \Einheit\hbox to0pt{\hskip\xcoord \Einheit
         \vrule height\Pfadd@cke width1 \Einheit depth\Pfadd@cke\hss}%
      \advance\xcoord by 1
    \else\ifnum#1=2
      \raise\ycoord \Einheit\hbox to0pt{\hskip\xcoord \Einheit
        \hbox{\hskip-\PfadD@cke\vrule height1 \Einheit width\PfadD@cke  
depth0pt}\hss}%
      \advance\ycoord by 1
    \else\ifnum#1=3
      \raise\ycoord \Einheit\hbox to0pt{\hskip\xcoord \Einheit
         \Line@(1,1){1}\hss}
      \advance\xcoord by 1
      \advance\ycoord by 1
    \else\ifnum#1=4
      \raise\ycoord \Einheit\hbox to0pt{\hskip\xcoord \Einheit
         \Line@(1,-1){1}\hss}
      \advance\xcoord by 1
      \advance\ycoord by -1
    \fi\fi\fi\fi
  \fi\next}
\def\hSSchritt{\leavevmode\raise-.4pt\hbox to0pt{\hss.\hss}\hskip.2\Einheit
  \raise-.4pt\hbox to0pt{\hss.\hss}\hskip.2\Einheit
  \raise-.4pt\hbox to0pt{\hss.\hss}\hskip.2\Einheit
  \raise-.4pt\hbox to0pt{\hss.\hss}\hskip.2\Einheit
  \raise-.4pt\hbox to0pt{\hss.\hss}\hskip.2\Einheit}
\def\vSSchritt{\vbox{\baselineskip.2\Einheit\lineskiplimit0pt
\hbox{.}\hbox{.}\hbox{.}\hbox{.}\hbox{.}}}
\def\DSSchritt{\leavevmode\raise-.4pt\hbox to0pt{%
  \hbox to0pt{\hss.\hss}\hskip.2\Einheit
  \raise.2\Einheit\hbox to0pt{\hss.\hss}\hskip.2\Einheit
  \raise.4\Einheit\hbox to0pt{\hss.\hss}\hskip.2\Einheit
  \raise.6\Einheit\hbox to0pt{\hss.\hss}\hskip.2\Einheit
  \raise.8\Einheit\hbox to0pt{\hss.\hss}\hss}}
\def\dSSchritt{\leavevmode\raise-.4pt\hbox to0pt{%
  \hbox to0pt{\hss.\hss}\hskip.2\Einheit
  \raise-.2\Einheit\hbox to0pt{\hss.\hss}\hskip.2\Einheit
  \raise-.4\Einheit\hbox to0pt{\hss.\hss}\hskip.2\Einheit
  \raise-.6\Einheit\hbox to0pt{\hss.\hss}\hskip.2\Einheit
  \raise-.8\Einheit\hbox to0pt{\hss.\hss}\hss}}
\def\SPfad(#1,#2),#3\endSPfad{\unskip\leavevmode
  \xcoord#1 \ycoord#2 \ZeichneSPfad#3\endSPfad}
\def\ZeichneSPfad#1{\ifx#1\endSPfad\let\next\relax
  \else\let\next\ZeichneSPfad
    \ifnum#1=1
      \raise\ycoord \Einheit\hbox to0pt{\hskip\xcoord \Einheit
         \hSSchritt\hss}%
      \advance\xcoord by 1
    \else\ifnum#1=2
      \raise\ycoord \Einheit\hbox to0pt{\hskip\xcoord \Einheit
        \hbox{\hskip-2pt \vSSchritt}\hss}%
      \advance\ycoord by 1
    \else\ifnum#1=3
      \raise\ycoord \Einheit\hbox to0pt{\hskip\xcoord \Einheit
         \DSSchritt\hss}
      \advance\xcoord by 1
      \advance\ycoord by 1
    \else\ifnum#1=4
      \raise\ycoord \Einheit\hbox to0pt{\hskip\xcoord \Einheit
         \dSSchritt\hss}
      \advance\xcoord by 1
      \advance\ycoord by -1
    \fi\fi\fi\fi
  \fi\next}
\def\Koordinatenachsen(#1,#2){\unskip
 \hbox to0pt{\hskip-.5pt\vrule height#2 \Einheit width.5pt depth1 \Einheit}%
 \hbox to0pt{\hskip-1 \Einheit \xcoord#1 \advance\xcoord by1
    \vrule height0.25pt width\xcoord \Einheit depth0.25pt\hss}}
\def\Koordinatenachsen(#1,#2)(#3,#4){\unskip
 \hbox to0pt{\hskip-.5pt \ycoord-#4 \advance\ycoord by1
    \vrule height#2 \Einheit width.5pt depth\ycoord \Einheit}%
 \hbox to0pt{\hskip-1 \Einheit \hskip#3\Einheit
    \xcoord#1 \advance\xcoord by1 \advance\xcoord by-#3
    \vrule height0.25pt width\xcoord \Einheit depth0.25pt\hss}}
\def\Gitter(#1,#2){\unskip \xcoord0 \ycoord0 \leavevmode
  \LOOP\ifnum\ycoord<#2
    \loop\ifnum\xcoord<#1
      \raise\ycoord \Einheit\hbox to0pt{\hskip\xcoord \Einheit\Punkt\hss}%
      \advance\xcoord by1
    \repeat
    \xcoord0
    \advance\ycoord by1
  \REPEAT}
\def\Gitter(#1,#2)(#3,#4){\unskip \xcoord#3 \ycoord#4 \leavevmode
  \LOOP\ifnum\ycoord<#2
    \loop\ifnum\xcoord<#1
      \raise\ycoord \Einheit\hbox to0pt{\hskip\xcoord \Einheit\Punkt\hss}%
      \advance\xcoord by1
    \repeat
    \xcoord#3
    \advance\ycoord by1
  \REPEAT}
\def\Label#1#2(#3,#4){\unskip \xdim#3 \Einheit \ydim#4 \Einheit
  \def\lo{\advance\xdim by-.5 \Einheit \advance\ydim by.5 \Einheit}%
  \def\llo{\advance\xdim by-.25cm \advance\ydim by.5 \Einheit}%
  \def\loo{\advance\xdim by-.5 \Einheit \advance\ydim by.25cm}%
  \def\o{\advance\ydim by.25cm}%
  \def\ro{\advance\xdim by.5 \Einheit \advance\ydim by.5 \Einheit}%
  \def\rro{\advance\xdim by.25cm \advance\ydim by.5 \Einheit}%
  \def\roo{\advance\xdim by.5 \Einheit \advance\ydim by.25cm}%
  \def\l{\advance\xdim by-.30cm}%
  \def\r{\advance\xdim by.30cm}%
  \def\lu{\advance\xdim by-.5 \Einheit \advance\ydim by-.6 \Einheit}%
  \def\llu{\advance\xdim by-.25cm \advance\ydim by-.6 \Einheit}%
  \def\luu{\advance\xdim by-.5 \Einheit \advance\ydim by-.30cm}%
  \def\u{\advance\ydim by-.30cm}%
  \def\ru{\advance\xdim by.5 \Einheit \advance\ydim by-.6 \Einheit}%
  \def\rru{\advance\xdim by.25cm \advance\ydim by-.6 \Einheit}%
  \def\ruu{\advance\xdim by.5 \Einheit \advance\ydim by-.30cm}%
  #1\raise\ydim\hbox to0pt{\hskip\xdim
     \vbox to0pt{\vss\hbox to0pt{\hss$#2$\hss}\vss}\hss}%
}
\catcode`\@=13

\def\ldreieck{\bsegment
  \rlvec(0.866025403784439 .5) \rlvec(0 -1)
  \rlvec(-0.866025403784439 .5)
  \savepos(0.866025403784439 -.5)(*ex *ey)
        \esegment
  \move(*ex *ey)
        }
\def\rdreieck{\bsegment
  \rlvec(0.866025403784439 -.5) \rlvec(-0.866025403784439 -.5)  \rlvec(0 1)
  \savepos(0 -1)(*ex *ey)
        \esegment
  \move(*ex *ey)
        }
\def\rhombus{\bsegment
  \rlvec(0.866025403784439 .5) \rlvec(0.866025403784439 -.5)
  \rlvec(-0.866025403784439 -.5)  \rlvec(0 1)
  \rmove(0 -1)  \rlvec(-0.866025403784439 .5)
  \savepos(0.866025403784439 -.5)(*ex *ey)
        \esegment
  \move(*ex *ey)
        }
\def\RhombusA{\bsegment
  \rlvec(0.866025403784439 .5) \rlvec(0.866025403784439 -.5)
  \rlvec(-0.866025403784439 -.5) \rlvec(-0.866025403784439 .5)
  \savepos(0.866025403784439 -.5)(*ex *ey)
        \esegment
  \move(*ex *ey)
        }
\def\RhombusB{\bsegment
  \rlvec(0.866025403784439 .5) \rlvec(0 -1)
  \rlvec(-0.866025403784439 -.5) \rlvec(0 1)
  \savepos(0 -1)(*ex *ey)
        \esegment
  \move(*ex *ey)
        }
\def\RhombusC{\bsegment
  \rlvec(0.866025403784439 -.5) \rlvec(0 -1)
  \rlvec(-0.866025403784439 .5) \rlvec(0 1)
  \savepos(0.866025403784439 -.5)(*ex *ey)
        \esegment
  \move(*ex *ey)
        }
\def\hdSchritt{\bsegment
  \lpatt(.05 .13)
  \rlvec(0.866025403784439 -.5)
  \savepos(0.866025403784439 -.5)(*ex *ey)
        \esegment
  \move(*ex *ey)
        }
\def\vdSchritt{\bsegment
  \lpatt(.05 .13)
  \rlvec(0 -1)
  \savepos(0 -1)(*ex *ey)
        \esegment
  \move(*ex *ey)
        }

\input texdraw

\begin{document}
\abovedisplayskip=6pt plus2pt minus 4pt
\belowdisplayskip=6pt plus2pt minus 4pt



\title[Enumeration of rhombus tilings]{\uppercase 
{{\bf Enumeration of rhombus tilings of a hexagon which contain a
fixed rhombus on its symmetry axis}}\\
{\rm\small (Extended Abstract)}}

\author{M.~Ciucu, M.~Fulmek and C.~Krattenthaler}
\address{
School of Mathematics\\
Institute for Advanced Study\\
Princeton, NJ 08540, USA\\
e-mail: ciucu@@ias.edu}
\address{
Institut f\"ur Mathematik der Universit\"at Wien\\
Strudlhofgasse 4, A-1090 Wien, Austria\\
e-mail: Mfulmek@@Mat.Univie.Ac.At, Kratt@@Pap.Univie.Ac.At
}



\maketitle
\def\abstractname{Summary}
\begin{abstract}
We compute the number
of rhombus tilings of a hexagon with sides $N,M,N,\break N,M,N$, which
contain a fixed rhombus on the symmetry axis. A special case
solves a problem posed by Jim Propp.
\end{abstract}


\section{Introduction}
In recent years, the enumeration of rhombus tilings of various regions
has attracted a lot of interest and was intensively studied, 
mainly because of the observation (see
\cite{KupeAA}) that the problem of enumerating all rhombus tilings of a hexagon
with sides $a,b,c,a,b,c$ and whose angles are $120^\circ$ 
(see Figure~\ref{fig:0}; throughout the paper by a
rhombus we always mean a rhombus with side lengths 1 and angles of
$60^\circ$ and $120^\circ$) 
is another way of stating the problem of
counting all plane partitions inside an $a\times b\times c$ box. The
latter problem was solved long ago by MacMahon 
\cite[Sec.~429, $q\rightarrow 1$; proof in Sec.~494]{MacMahon}. Therefore:

\smallskip
{\it The number of all rhombus tilings of a hexagon
with sides $a,b,c,a,b,c$ equals}
\begin{equation}
\label{eq:MacMahon}
\prod_{i=1}^a\prod_{j=1}^b\prod_{k=1}^c\frac{i+j+k-1}{i+j+k-2}.
\end{equation}
(The form of the expression is due to Macdonald.)
\smallskip

\begin{figure}
\centertexdraw{
  \drawdim truecm  \linewd.02
  \rhombus \rhombus \rhombus \rhombus \ldreieck
  \move (-0.866025403784439 -.5)
  \rhombus \rhombus \rhombus \rhombus \rhombus \ldreieck
  \move (-1.732050807568877 -1)
  \rhombus \rhombus \rhombus \rhombus \rhombus \rhombus \ldreieck
  \move (-1.732050807568877 -1)
  \rdreieck
  \rhombus \rhombus \rhombus \rhombus \rhombus \rhombus \ldreieck
  \move (-1.732050807568877 -2)
  \rdreieck
  \rhombus \rhombus \rhombus \rhombus \rhombus \rhombus \ldreieck
  \move (-1.732050807568877 -3)
  \rdreieck
  \rhombus \rhombus \rhombus \rhombus \rhombus \rhombus 
  \move (-1.732050807568877 -4)
  \rdreieck
  \rhombus \rhombus \rhombus \rhombus \rhombus 
  \move (-1.732050807568877 -5)
  \rdreieck
  \rhombus \rhombus \rhombus \rhombus 
\move(8 0)
\bsegment
  \drawdim truecm  \linewd.02
  \rhombus \rhombus \rhombus \rhombus \ldreieck
  \move (-0.866025403784439 -.5)
  \rhombus \rhombus \rhombus \rhombus \rhombus \ldreieck
  \move (-1.732050807568877 -1)
  \rhombus \rhombus \rhombus \rhombus \rhombus \rhombus \ldreieck
  \move (-1.732050807568877 -1)
  \rdreieck
  \rhombus \rhombus \rhombus \rhombus \rhombus \rhombus \ldreieck
  \move (-1.732050807568877 -2)
  \rdreieck
  \rhombus \rhombus \rhombus \rhombus \rhombus \rhombus \ldreieck
  \move (-1.732050807568877 -3)
  \rdreieck
  \rhombus \rhombus \rhombus \rhombus \rhombus \rhombus 
  \move (-1.732050807568877 -4)
  \rdreieck
  \rhombus \rhombus \rhombus \rhombus \rhombus 
  \move (-1.732050807568877 -5)
  \rdreieck
  \rhombus \rhombus \rhombus \rhombus 
  \linewd.08
  \move(0 0)
  \RhombusA \RhombusB \RhombusB 
  \RhombusA \RhombusA \RhombusB \RhombusA \RhombusB \RhombusB
  \move (-0.866025403784439 -.5)
  \RhombusA \RhombusB \RhombusB \RhombusB \RhombusB
  \RhombusA \RhombusA \RhombusB \RhombusA 
  \move (-1.732050807568877 -1)
  \RhombusB \RhombusB \RhombusA \RhombusB \RhombusB \RhombusA
  \RhombusB \RhombusA \RhombusA 
  \move (1.732050807568877 0)
  \RhombusC \RhombusC \RhombusC 
  \move (1.732050807568877 -1)
  \RhombusC \RhombusC \RhombusC 
  \move (3.464101615137755 -3)
  \RhombusC 
  \move (-0.866025403784439 -.5)
  \RhombusC
  \move (-0.866025403784439 -1.5)
  \RhombusC
  \move (0.866025403784439 -2.5)
  \RhombusC \RhombusC 
  \move (0.866025403784439 -3.5)
  \RhombusC \RhombusC \RhombusC 
  \move (2.598076211353316 -5.5)
  \RhombusC 
  \move (0.866025403784439 -5.5)
  \RhombusC 
  \move (-1.732050807568877 -3)
  \RhombusC 
  \move (-1.732050807568877 -4)
  \RhombusC 
  \move (-1.732050807568877 -5)
  \RhombusC \RhombusC 
\esegment
\htext (-1.5 -9){\small a. A hexagon with sides $a,b,c,a,b,c$,}
\htext (-1.5 -9.5){\small \hphantom{a. }where $a=3$, $b=4$, $c=5$}
\htext (6.8 -9){\small b. A rhombus tiling of a hexagon}
\htext (6.8 -9.5){\small \hphantom{b. }with sides $a,b,c,a,b,c$}
\rtext td:0 (4.3 -4){$\sideset {}  c 
    {\left.\vbox{\vskip2.2cm}\right\}}$}
\rtext td:60 (2.6 -.5){$\sideset {} {} 
    {\left.\vbox{\vskip1.7cm}\right\}}$}
\rtext td:120 (-.44 -.25){$\sideset {}  {}  
    {\left.\vbox{\vskip1.3cm}\right\}}$}
\rtext td:0 (-2.3 -3.6){$\sideset {c}  {} 
    {\left\{\vbox{\vskip2.2cm}\right.}$}
\rtext td:240 (0 -7){$\sideset {}  {}  
    {\left.\vbox{\vskip1.7cm}\right\}}$}
\rtext td:300 (3.03 -7.25){$\sideset {}  {}  
    {\left.\vbox{\vskip1.4cm}\right\}}$}
\htext (-.9 0.1){$a$}
\htext (2.8 -.1){$b$}
\htext (3.2 -7.8){$a$}
\htext (-0.3 -7.65){$b$}
}
\caption{}
\label{fig:0}
\end{figure}

A statistical investigation of which rhombi lie in a random rhombus 
tiling has been undertaken, on an {\it asymptotic} level, by Cohn, 
Larsen and Propp \cite{CoLPAA}. On the {\it exact\/} (enumerative) level, 
Propp \cite[Problem~1]{PropAA} observed numerically that apparently exactly one
third of 
the rhombus tilings of a hexagon with side lengths 
$2n-1,2n,2n-1,2n-1,2n,2n-1$ contain the central rhombus.

In this article we present the solution of an even more general problem,
namely the enumeration of all rhombus tilings of a hexagon with side lengths
$N,M,N,N,M,N$ which contain an {\it arbitrary} fixed rhombus
on the symmetry axis which cuts through the sides of length $M$ (see
Figure~\ref{fig:1} for illustration; the fixed rhombus is shaded). Our
results are the following.

\begin{figure}
\bigskip
\centertexdraw{
  \drawdim truecm  \linewd.02
  \rhombus \rhombus \rhombus \ldreieck
  \move (-0.866025403784439 -.5)
  \rhombus \rhombus \rhombus \rhombus \ldreieck
  \move (-1.732050807568877 -1)
  \rhombus \rhombus \rhombus \rhombus \rhombus
  \move (-1.732050807568877 -1)
  \rdreieck
  \rhombus \rhombus \rhombus \rhombus
  \move (-1.732050807568877 -2)
  \rdreieck
  \rhombus \rhombus \rhombus
  \move (-0.866025403784439 -1.5)
  \bsegment
    \rlvec(0.866025403784439 -.5) \rlvec(-0.866025403784439 -.5)
	   \rlvec(-0.866025403784439 .5)
    \lfill f:0.8
  \esegment
  \linewd.05
  \move (-1.732050807568877 -2)
  \rlvec(5.196152422706631 0)
\move(8 0)
\bsegment
  \drawdim truecm  \linewd.02
  \rhombus \rhombus \rhombus \ldreieck
  \move (-0.866025403784439 -.5)
  \rhombus \rhombus \rhombus \rhombus \ldreieck
  \move (-1.732050807568877 -1)
  \rhombus \rhombus \rhombus \rhombus \rhombus \ldreieck
  \move (-1.732050807568877 -1)
  \rdreieck \rhombus \rhombus \rhombus \rhombus \rhombus
  \move (-1.732050807568877 -2)
  \rdreieck \rhombus \rhombus \rhombus \rhombus
  \move (-1.732050807568877 -3)
  \rdreieck \rhombus \rhombus \rhombus
  \move (1.732050807568877 -2)
  \bsegment
    \rlvec(0.866025403784439 -.5) \rlvec(-0.866025403784439 -.5)
	   \rlvec(-0.866025403784439 .5)
    \lfill f:0.8
  \esegment
  \linewd.05
  \move (-1.732050807568877 -2.5)
  \rlvec(5.196152422706631 0)
\esegment
\htext (-2 -6){\small A hexagon with sides $N,M,N,$}
\htext (-2 -6.5){\small $N,M,N$ and fixed rhombus }
\htext (-2 -7){\small $l$, where $N=3$, $M=2$, $l=1$.}
\htext (-2 -7.5){\small The thick  horizontal line indicates}
\htext (-2 -8){\small the symmetry axis.}
\htext (5.8 -6){\small A hexagon with sides $N,M,N,$}
\htext (5.8 -6.5){\small $N,M,N$ and fixed rhombus }
\htext (5.8 -7){\small $l$, where $N=3$, $M=3$, $l=2$.}
\htext (5.8 -7.5){\small The thick  horizontal line indicates}
\htext (5.8 -8){\small the symmetry axis.}
\htext (-.9 0.1){$N$}
\htext (2.4 0.1){$N$}
\htext (-2.5 -2){$M$}
\htext (7.1 0.1){$N$}
\htext (10.4 0.1){$N$}
\htext (5.5 -2.6){$M$}
}
\caption{}
\label{fig:1}
\end{figure}

\begin{thm}
\label{thm:MEven}
Let $m$ be a nonnegative integer and $N$ be a positive integer. The number
of rhombus tilings of a hexagon with sides $N,2m,N,N,2m,N$, which
contain the $l$-th rhombus on the symmetry axis 
which cuts through the sides of length $2m$, equals
\begin{multline}
\label{eq:MEven}
\frac{m\binom{m+N}{m}\binom{m+N-1}{m}}{\binom{2m+2N-1}{2m}}
\sum_{e=0}^{l-1}(-1)^{e}\binom{N}{e}
	\frac{
		(N-2e)\po{\frac{1}{2}}{e}
	}{
		(m+e)(m+N-e)\po{\frac{1}{2}-N}{e}
	}\\
\times
\prod_{i=1}^N\prod_{j=1}^N\prod_{k=1}^{2m}\frac{i+j+k-1}{i+j+k-2},
\end{multline}
where the shifted factorial
$(a)_k$ is given by $(a)_k:=a(a+1)\cdots(a+k-1)$,
$k\ge1$, $(a)_0:=1$.
\end{thm}

\begin{thm}
\label{thm:MOdd}
Let $m$ and $N$ be positive integers. The number
of rhombus tilings of a hexagon with sides $N+1,2m-1,N+1,N+1,2m-1,N+1$, which
contain the $l$-th rhombus on the symmetry axis
which cuts through the sides of length $2m-1$, equals
\begin{multline}
\label{eq:MOdd}
\frac{m\binom{m+N}{m}\binom{m+N-1}{m}}{\binom{2m+2N-1}{2m}}
\sum_{e=0}^{l-1}(-1)^{e}\binom{N}{e}
	\frac{
		(N-2e)\po{\frac{1}{2}}{e}
	}{
		(m+e)(m+N-e)\po{\frac{1}{2}-N}{e}
	}\\
\times
\prod_{i=1}^{N+1}\prod_{j=1}^{N+1}\prod_{k=1}^{2m-1}\frac{i+j+k-1}{i+j+k-2}.
\end{multline}
\end{thm}

The special case of Theorem~\ref{thm:MEven} where the fixed rhombus is the 
central rhombus was proved by the first and third
author \cite{CiucKratAB},
and independently by Helfgott and Gessel \cite[Theorem~2]{HeGeAA},
using a different method.
Building on the approach of \cite{CiucKratAB}, the second and third author 
\cite{FuKrAC}
were able to generalize the enumeration to the above theorems.

The special case $N=2n-1$, $m=n$ of Theorem~\ref{thm:MEven} does indeed
imply Propp's conjecture.

\begin{cor}
\label{cor:ProppVerm}
Let $n$ be a positive integer.
Exactly one third of the rhombus tilings of a
hexagon with sides $2n-1,2n,2n-1,2n-1,2n,2n-1$  cover the
central rhombus. The same is true for a hexagon with
sides $2n,2n-1,2n,2n,2n-1,2n$.
\end{cor}

Finally, from Theorems~\ref{thm:MEven} and \ref{thm:MOdd}, we derive
an ``arcsine law" for this kind of enumeration. It complements the 
asymptotic results by Cohn, Larsen and Propp \cite{CoLPAA}.

\begin{thm}\label{thm:Asymptotics}
Let $a$ be any nonnegative real number, let $b$ be a real number with 
$0<b<1$.
For $m\sim a N$ and $l\sim b N$, the proportion of rhombus tilings of a
hexagon with sides $N,2m,N,N,2m,N$ or $N+1,2m-1,N+1,N+1,2m-1,N+1$, 
which contain the
$l$-th rhombus on the symmetry axis which cuts through the sides of length
$2m$, respectively $2m-1$, in the total
number of rhombus tilings is asymptotically
\begin{equation}
\label{eq:Asymptotics}
\frac{2}{\pi}\arcsin\left(
	\frac{\sqrt{b(1-b)}}{\sqrt{(a+b)(a-b+1)}}
\right)
\end{equation}
as $N$ tends to infinity.
\end{thm}

In the remainder of this article we sketch proofs of these results.
In the next section we provide brief outlines of the proofs of 
Theorems~\ref{thm:MEven}, \ref{thm:MOdd}, \ref{thm:Asymptotics} and
Corollary~\ref{cor:ProppVerm}. The 
proof (or rather, a sketch of the proof) 
of a crucial auxiliary lemma is deferred to Section~3.

\section{Outline of proofs}

{\sc Outline of proof of Theorems~\ref{thm:MEven} and \ref{thm:MOdd}}.
The proofs of both Theorems are very similar. We will mainly concentrate
on the proof of Theorem~\ref{thm:MEven}.

There are four basic steps. 

\begin{figure}
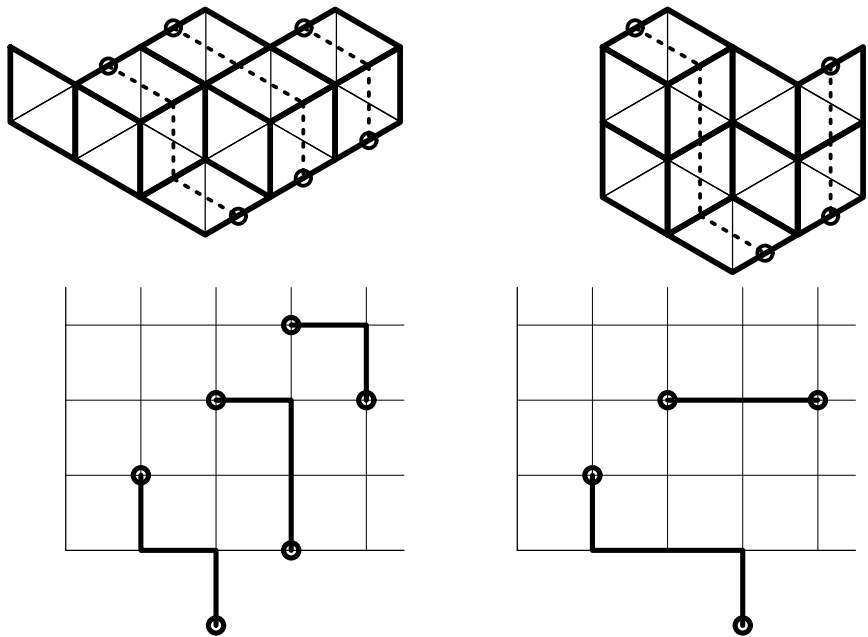

\bigskip
\centertexdraw{
\drawdim truecm  \linewd.02
\move(0 -1)
\bsegment
  \rhombus \rhombus \rhombus \ldreieck
  \move (-0.866025403784439 -.5)
  \rhombus \rhombus \rhombus
  \move (-1.732050807568877 -1)
  \rhombus \rhombus
  \move (-1.732050807568877 -1)
  \rdreieck
  \move (1.732050807568877 -4)
  \rhombus \ldreieck
  \move (-0 -4)
  \rhombus \rhombus \rhombus
  \move (-.866025403784439 -4.5)
  \rhombus \rhombus \rhombus
  \move (-1.732050807568877 -4)
  \rdreieck \rhombus \rhombus \rhombus
  \move (-1.732050807568877 -1)
  \bsegment
    \rlvec(0.866025403784439 .5) \rlvec(0 -1)
	   \rlvec(-0.866025403784439 -.5)
    \lfill f:0.5
  \esegment
  \move (2.598076211353316 -.5)
  \bsegment
    \rlvec(0.866025403784439 -.5) \rlvec(0 -1)
	   \rlvec(-0.866025403784439 .5)
    \lfill f:0.5
  \esegment
  \linewd.08
  \move (-1.732050807568877 -1)
  \RhombusB
  \move (2.598076211353316 -.5)
  \RhombusC
  \htext (-2.0 0.6){
  	$\overbrace{\hphantom{xxxxxxxxxxxxxxxxxxxxxxx}}^{S^\prime(3,1)}$
	}
  \htext (-0.8 -2.7){$\underbrace{\hphantom{xxxxxxxxxxxxxx}}_{S(2,1)}$}
  \htext (-2.0 -7.2){
  	$\underbrace{\hphantom{xxxxxxxxxxxxxxxxxxxxxxx}}_{C(3,1,1)}$
	}
\esegment
\move(7 -1)
\bsegment
  \drawdim truecm  \linewd.02
  \rhombus \rhombus \rhombus \ldreieck
  \move (-0.866025403784439 -.5)
  \rhombus \rhombus \rhombus \rhombus
  \move (-1.732050807568877 -1)
  \rhombus \rhombus \rhombus
  \move (-1.732050807568877 -1)
  \rdreieck \rhombus
  \move (2.598076211353316 -4.5)
  \ldreieck
  \move (1.732050807568877 -5)
  \rhombus \ldreieck
  \move (-0.866025403784439 -4.5)
  \rhombus \rhombus \rhombus \rhombus
  \move (-1.732050807568877 -4)
  \rdreieck \rhombus \rhombus \rhombus \rhombus
  \move (-1.732050807568877 -5)
  \rdreieck \rhombus \rhombus \rhombus
  \move (-1.732050807568877 -4)
  \bsegment
    \rlvec(0.866025403784439 -.5) \rlvec(0 -2)
	   \rlvec(-0.866025403784439 .5)
    \lfill f:0.5
  \esegment
  \move (3.464101615137754 -4)
  \bsegment
    \rlvec(-0.866025403784439 -.5) \rlvec(0 -2)
	   \rlvec(0.866025403784439 .5)
    \lfill f:0.5
  \esegment
  \linewd.08
  \move (-1.732050807568877 -4)
  \RhombusC
  \move (-1.732050807568877 -5)
  \RhombusC
  \move (2.598076211353316 -4.5)
  \RhombusB
  \move (2.598076211353316 -5.5)
  \RhombusB
  \htext (-2.0 0.6){
  	$\overbrace{\hphantom{xxxxxxxxxxxxxxxxxxxxxxx}}^{S(3,1)}$
	}
  \htext (-2.0 -3.8){
  	$\overbrace{\hphantom{xxxxxxxxxxxxxxxxxxxxxxx}}^{C^\prime(3,1,1)}$
	}
  \htext (-0.9 -8.2){
  	$\underbrace{\hphantom{xxxxxxxxxxxxxx}}_{C(2,1,1)}$
	}
\esegment
\move (-1.5 -9.7)
\bsegment
	\htext (0 0){\small The hexagons from Figure~\ref{fig:1} cut
along the symmetry axis, according}
	\htext (0 -.5){\small to the Matchings Factorization Theorem. In the
shaded regions,}
	\htext (0 -1){\small ``forced" rhombi are shown with thick lines.}
\esegment
}
\caption{Hexagons, cut in two}
\label{fig:2}
\end{figure}

{\it Step~1. Application of the Matchings Factorization Theorem}.
First, rhombus tilings of the hexagon with sides $N,2m,N,N,2m,N$
can be interpreted as perfect matchings of the dual graph of the
triangulated hexagon, i.e., the (bipartite) graph $G(V,E)$, where
the set of vertices $V$ consists of the triangles of the hexagon's
triangulation, and where
two vertices are connected by an edge if the corresponding triangles
are adjacent. Enumerating only those rhombus tilings which contain
a fixed rhombus, under this translation amounts to enumerating only
those perfect matchings which contain the edge corresponding to this
rhombus, or, equivalently, we may consider just perfect matchings of the
graph which results from $G(V,E)$ by removing this edge. 
Clearly, since the fixed rhombus was located on the symmetry axis,
this graph is symmetric.
Hence, we may apply
the first author's Matchings Factorization
Theorem \cite[Thm.~1.2]{CiucAB}. In general, this theorem says that 
the number of perfect matchings of a symmetric graph $G$ equals a certain
power of 2 times the number of perfect matchings of a graph $G^+$ (which
is, roughly speaking, the ``upper half" of $G$) times a weighted count 
of perfect
matchings of a graph $G^-$ (which is, roughly speaking, the ``lower half" of
$G$), in which the edges on the symmetry axis count with weight $1/2$ only.
Applied to our case, and retranslated into rhombus tilings, the Matchings
Factorization Theorem implies the following:

{\em The number
of rhombus tilings of a hexagon with sides $N,2m,N,N,2m,N$, which
contain the $l$-th rhombus on the symmetry axis 
which cuts through the sides of length $2m$, equals
\begin{equation}\label{eq:MatchFact}
2^{N-1}R(S'(N,m))\tilde R(C(N,m,l)),
\end{equation}}
where $S'(N,m)$ denotes the ``upper half" of our hexagon with the fixed
rhombus removed (see Figure~\ref{fig:2}), where $R(S'(m,n))$ denotes 
the number of rhombus tilings of $S'(m,n)$, where $C(N,m,l)$ denotes the 
``lower half" (again, see Figure~\ref{fig:2}), and where $\tilde R(C(m,n,l))$
denotes the weighted count of rhombus tilings of $C(m,n,l)$ in which
each of the top-most (horizontal) rhombi counts with weight $1/2$.
(Both, $S'(N,m)$ and $C(N,m,l)$ are roughly pentagonal. The notations $S'(N,m)$
and $C(N,m,l)$ stand for ``simple part" and ``complicated part", respectively,
as it will turn out that the count $R(S'(N,m))$ will be rather 
straight-forward, while the count $\tilde R(C(N,m,l))$ will turn out be
considerably harder.)

In the case of Theorem~\ref{thm:MEven} it 
is immediately obvious, that the rhombi along the left-most and
right-most vertical strip of $S'(N,m)$ must be contained in any rhombus 
tiling of $S'(N,m)$. Hence, we may safely remove these strips. Let us denote
the resulting region by $S(N-1,m)$. From \eqref{eq:MatchFact}
we obtain that
{\em the number
of rhombus tilings of a hexagon with sides $N,2m,N,N,2m,N$, which
contain the $l$-th rhombus on the symmetry axis 
which cuts through the sides of length $2m$, equals
\begin{equation}
\label{eq:CiucuMEven}
2^{N-1}R(S(N-1,m))\tilde R(C(N,m,l)).
\end{equation}}

Similarly, for the case of Theorem~\ref{thm:MOdd}, we obtain that
{\em the number
of rhombus tilings of a hexagon with sides $N,2m-1,N,N,2m-1,N$, which
contain the $l$-th rhombus on the symmetry axis 
which cuts through the sides of length $2m-1$, equals
\begin{equation}
\label{eq:CiucuMOdd}
2^{N-1}R(S(N,m-1))\tilde R(C(N-1,m,l)).
\end{equation}}

\smallskip
{\it Step~2. From rhombus tilings to nonintersecting lattice paths.}
There is a standard translation from rhombus tilings to nonintersecting
lattice paths. We apply it to our regions $S(N,m)$ and $C(N,m,l)$. 
Figure~\ref{fig:3} illustrates this translation for the (``complicated")
lower parts in Figure~\ref{fig:2}.

\begin{figure}
\bigskip
\centertexdraw{
  \drawdim truecm  \linewd.02
\move(0 0)
\bsegment
\bsegment
  \move (1.732050807568877 0)
  \rhombus \ldreieck
  \move (-0 0)
  \rhombus \rhombus \rhombus
  \move (-.866025403784439 -0.5)
  \rhombus \rhombus \rhombus
  \move (-1.732050807568877 -0)
  \rdreieck \rhombus \rhombus \rhombus
  \linewd.08
  \move (-1.732050807568877 0)
  \RhombusC \RhombusC
  \rmove(-0.866025403784439 0.5)
  \RhombusA \RhombusB \RhombusA
  \rmove(-0.866025403784439 2.5)
  \RhombusA \RhombusA \RhombusB
  \rmove (0 2)
  \RhombusA \RhombusB
  \linewd.05
  \move (-0.4330127018922194 -.25)
  \lcir r:.1
  \hdSchritt \vdSchritt \hdSchritt
  \lcir r:.1
  \move (0.4330127018922194 .25)
  \lcir r:.1
  \hdSchritt \hdSchritt \vdSchritt
  \lcir r:.1
  \move (2.165063509461096 .25)
  \lcir r:.1
  \hdSchritt \vdSchritt
  \lcir r:.1
\esegment
\move(7 0)
\bsegment
  \linewd.02
  \move (1.732050807568877 -0.5)
  \ldreieck
  \move (-0.866025403784439 -0)
  \rhombus \rhombus \rhombus \ldreieck
  \move (-0.866025403784439 -0)
  \rdreieck \rhombus \rhombus \rhombus
  \move (-0.866025403784439 -1)
  \rdreieck \rhombus \rhombus
  \linewd.08
  \move (-0.866025403784439 0)
  \RhombusC
  \rmove (-0.866025403784439 -0.5)
  \RhombusC
  \move (-0.866025403784439 0)
  \RhombusA \RhombusB \RhombusB \RhombusA
  \move (0.866025403784439 0)
  \RhombusC
  \rmove (-0.866025403784439 -0.5)
  \RhombusC
  \move (1.732050807568877 -0.5)
  \RhombusB \RhombusB
  \linewd.05
  \move (-0.4330127018922194 .25)
  \lcir r:.1
  \hdSchritt \vdSchritt \vdSchritt \hdSchritt
  \lcir r:.1
  \move (2.165063509461096 -.25)
  \lcir r:.1
  \vdSchritt \vdSchritt
  \lcir r:.1
\esegment
\esegment
\move(0 -6.7)
\bsegment
	\move(-1 0)
	\bsegment
		\linewd.02
		\move (0 0) \lvec (4.5 0)
		\move (0 0) \lvec (0 3.5)
		\linewd.01
		\move (0 1) \rlvec (4.5 0)
		\move (0 2) \rlvec (4.5 0)
		\move (0 3) \rlvec (4.5 0)
		\move (1 0) \rlvec (0 3.5)
		\move (2 0) \rlvec (0 3.5)
		\move (3 0) \rlvec (0 3.5)
		\move (4 0) \rlvec (0 3.5)
		\linewd.07
		\move (1 1) \lcir r:.1 \rlvec(0 -1) \rlvec(1 0) \rlvec(0 -1)  
\lcir r:.1
		\move (2 2) \lcir r:.1 \rlvec(1 0) \rlvec(0 -1) \rlvec(0 -1)  
\lcir r:.1
		\move (3 3) \lcir r:.1 \rlvec(1 0) \rlvec(0 -1) \lcir r:.1
	\esegment
	\move(5 0)
	\bsegment
		\linewd.02
		\move (0 0) \lvec (4.5 0)
		\move (0 0) \lvec (0 3.5)
		\linewd.01
		\move (0 1) \rlvec (4.5 0)
		\move (0 2) \rlvec (4.5 0)
		\move (0 3) \rlvec (4.5 0)
		\move (1 0) \rlvec (0 3.5)
		\move (2 0) \rlvec (0 3.5)
		\move (3 0) \rlvec (0 3.5)
		\move (4 0) \rlvec (0 3.5)
		\linewd.07
		\move (1 1) \lcir r:.1 \rlvec(0 -1) \rlvec(1 0) \rlvec(1 0)
				\rlvec(0 -1) \lcir r:.1
		\move (2 2) \lcir r:.1 \rlvec(1 0) \rlvec(1 0) \lcir r:.1
	\esegment
\esegment
\move(-3.5 -8.5)
\bsegment
	\htext (0 0){\small Tilings for the ``complicated parts'' from
		Figure~\ref{fig:2}, interpreted as lattice paths:}
	\htext (0 -.5){\small Reflect and rotate the dotted paths to obtain
		the paths in the lower
		picture.}
\esegment
}
\caption{Lattice path interpretation}
\label{fig:3}
\end{figure}

For the ``simple" pentagonal part $S(N,m)$ we obtain the following:
The number\break $R(S(N,m))$ 
of rhombus tilings of $S(N,m)$ equals the number of families
$(P_1,P_2,\dots,P_N)$ of nonintersecting lattice paths consisting of
horizontal unit steps in the positive direction and vertical unit steps
in the negative direction, where $P_i$ runs from $(i,i)$ to $(i+m,2i-N-1)$,
$i=1,2,\dots,N$.

Similarly, for the ``complicated" pentagonal part $C(N,m,l)$ we obtain:
The weighted count $\tilde R(C(N,m,l))$ 
of rhombus tilings of $C(N,m,l)$ equals the weighted count of families
$(P_1,P_2,\dots,P_N)$ of nonintersecting lattice paths consisting of
horizontal unit steps in the positive direction and vertical unit steps
in the negative direction, where $P_i$ runs from 
$(i,i)$ to $(i+m,2i-N-1)$ if $i\neq l$, while
$P_l$ runs from $(l,l)$ to $(l+m,2l-N)$; with the additional twist
that path $P_i$ ($i\neq l$) has weight $1/2$ if it ends with a
vertical step.

\smallskip
{\it Step~3. From nonintersecting lattice paths to determinants}.
Now, by using the main theorem on nonintersecting lattice paths 
\cite[Cor.~2]{GeViAB} (see also \cite[Theorem~1.2]{StemAE}),
we may write $R(S(N,m))$ and $\tilde R(C(N,m,l))$ as
determinants. Namely, 
we have
\begin{align}
R(S(N,m))	
&=\det_{1\leq i,j\leq N}\left(\binom{N+m-i+1}{m+i-j}\right)\label{eq:SNm1},
\end{align}
and
\begin{equation}
\label{eq:CNml1}
\tilde R(C(N,m,l))	=
	\det_{1\leq i,j\leq N}\left(
		\begin{cases}
\frac{(N+m-i)!}{(m+i-j)!\,(N+j-2i+1)!}(m+\frac{N-j+1}{2})&\text{ if }i\neq l\\
			\frac{(N+m-i)!}{(m+i-j)!\,(N+j-2i)!}&\text{ if }i=l
		\end{cases}
	\right).
\end{equation}

\smallskip
{\it Step~4. Determinant evaluations}. Clearly, once we are able
to evaluate the determinants in \eqref{eq:SNm1} and \eqref{eq:CNml1},
Theorems~\ref{thm:MEven} and \ref{thm:MOdd} will immediately follow
from \eqref{eq:CiucuMEven} and \eqref{eq:CiucuMOdd}, respectively,
upon routine simplification.
Indeed, for the determinant in \eqref{eq:SNm1} we have the following.
\begin{lem}
\label{lem:simplepart}
\begin{align}
&\det_{1\leq i,j\leq N}\left(\binom{N+m-i+1}{m+i-j}\right)
 =\prod_{i=1}^{N}\frac{(N+m-i+1)!\,(i-1)!\,\po{2m+i+1}{i-1}}{(m+i-1)!\,(2N-2i+1)!}.
			\label{eq:SNm2}
\end{align}
\end{lem}
\begin{pf}This determinant evaluation follows without difficulty from
a determinant lemma in \cite[Lemma~2.2]{KratAM}. The corresponding
computation is contained in the proof of Theorem~5 in \cite{KratAK}
(set $r=N$, $\lambda_s=m$, $B=2$, $a+\alpha-b=2m$ there, and then reverse
the order of rows and columns). 
\end{pf}

On the other hand, the determinant in \eqref{eq:CNml1} evaluates as follows.
\begin{lem}
\label{lem:complexpart}
\begin{multline}
\label{eq:CNml2}
	\det_{1\leq i,j\leq N}\left(
		\begin{cases}
\frac{(N+m-i)!}{(m+i-j)!\,(N+j-2i+1)!}(m+\frac{N-j+1}{2})&\text{ if }i\neq l\\
			\frac{(N+m-i)!}{(m+i-j)!\,(N+j-2i)!}&\text{ if }i=l
		\end{cases}
	\right)\\
= \prod_{i=1}^{N}\frac{(N+m-i)!}{(m+i-1)!\,(2N-2i+1)!}
\prod _{i=1}^{\floor{N/2}}
	\left(
		\po{m+i}{N-2i+1}\,\po{m+i+\frac{1}{2}}{N-2i}
	\right)\kern1.5cm
\\
\times
2^{\frac{(N-1)(N-2)}{2}}\frac{
	\po{m}{N+1}\prod_{j=1}^{N}(2j-1)!
}{
	N!\prod_{i=1}^{\floor{\frac{N}{2}}}\po{2i}{2n-4i+1}
}
\sum_{e=0}^{l-1}
(-1)^{e}\binom{N}{e}\frac{
	(N-2e)\,\po{\frac{1}{2}}{e}
}{
	(m+e)\,(m+N-e)\,\po{\frac{1}{2}-N}{e}
}.
\end{multline}
\end{lem}
This determinant evaluation is much more complex than the determinant
evaluation of Lemma~\ref{lem:simplepart}, and, as such, is the most
difficult part in our derivation of Theorems~\ref{thm:MEven} and
\ref{thm:MOdd}. We provide a sketch of how to evaluate this determinant
in the next section.

\smallskip
Altogether, Steps~1--4 establish Theorems~\ref{thm:MEven} and
\ref{thm:MOdd}.
\hfil\qed

\bigskip
{\sc Proof of Corollary~\ref{cor:ProppVerm}}.
We have to compute the ratio of the expression \eqref{eq:MEven}, 
with $N=2n-1$, $m=n$, by
the expression \eqref{eq:MacMahon}, with $a=b=2n-1$, $c=2n$, respectively 
the ratio of the expression \eqref{eq:MOdd}, with $N=2n-1$, $m=n$, by
the expression \eqref{eq:MacMahon}, with $a=b=2n$, $c=2n-1$. 
Clearly, except for
trivial manipulations, we will be done once we are able to evaluate
the sum in \eqref{eq:MEven} (which is the same as the one in 
\eqref{eq:MOdd}) for $N=2n-1$, $m=n$, and $l=n$.

We claim that 
$$\sum_{e=0}^{n-1}(-1)^{e}\binom{2n-1}{e}
	\frac{
		(2n-2e-1)\po{\frac{1}{2}}{e}
	}{
		(n+e)(3n-e-1)\po{\frac{3}{2}-2n}{e}
	}=2^{n-1}\frac {n!\,(n-1)!\,(6n-3)!!} {(3n)!\,(4n-3)!!}.
$$
Let us denote the sum by $S(n)$.
Then an application of the Gosper--Zeilberger algorithm 
\cite{PeWZAA,ZeilAM,ZeilAV} (we used the {\sl Mathematica} implementation
by Paule and Schorn \cite{PaScAA}) yields the relation
$$2n(2n+1)(6n-1)(6n+1)S(n)-(3n+1)(3n+2)(4n-1)(4n+1)S(n+1)=0,$$
which easily proves the claimed summation by an induction on $n$.
\hfill\qed

\bigskip
{\sc Outline of proof of Theorem~\ref{thm:Asymptotics}}.
{}From MacMahon's formula~\eqref{eq:MacMahon} for the total number of rhombus
tilings together with Theorems~\ref{thm:MEven} and \ref{thm:MOdd} we deduce
immediately that the proportion is indeed the same for both cases
$N,2m,N$ and $N+1,2m-1,N+1$, and that it is given by
\begin{equation}
\label{eq:proportion}
\frac{m\binom{m+N}{m}\binom{m+N-1}{m}}{\binom{2m+2N-1}{2m}}
\sum_{e=0}^{l-1}(-1)^{e}\binom{N}{e}
	\frac{
		(N-2e)\,\po{\frac{1}{2}}{e}
	}{
		(m+e)\,(m+N-e)\,\po{\frac{1}{2}-N}{e}
	}.
\end{equation}

We write the sum in \eqref{eq:proportion} in a hypergeometric fashion,
to get
\begin{equation}
\label{eq:proportion-hyp}
\frac{(2N-1)!\,\left(\po{m+1}{N-1}\right)^2}{(N-1)!^2\,\po{2m+1}{2N-1}}
\sum _{e=0} ^{l-1}\frac {(-N)_e\,(1-\frac{N}{2})_e\,
(m)_e\,(-m-N)_e\,(\frac{1}{2})_e} 
{(-\frac{N}{2})_e\,(1-m-N)_e\,(1+m)_e\,(\frac{1}{2}-N)_e\,e!}.
\end{equation}
Using a special case of Whipple's transformation (see \cite[(2.4.1.1)]{SlatAC}),
we transform the sum in \eqref{eq:proportion-hyp} into a 
$_4F_3$-series, thus obtaining
\begin{equation}\notag
\frac{(2N-1)!\,\left(\po{m+1}{N-1}\right)^2}{(N-1)!^2\,\po{2m+1}{2N-1}}
\frac{\po{-N+1}{l-1}\,\po{-l+\frac{1}{2}}{l-1}}%
{\po{-N+\frac{1}{2}}{l-1}\,\po{-l+1}{l-1}}
\HypF{4}{3}{1,\frac{1}{2},l-N,1-l}{1+m,1-m-N,\frac{3}{2}}{1}
\end{equation}
for the ratio \eqref{eq:proportion}.

Next we apply Bailey's transformation between
two balanced $\HypFsimple{4}{3}$-series (see\break \cite[(4.3.5.1)]{SlatAC}),
which gives the expression
\begin{multline}\label{eq:fact4F3}
\frac{(2l)!\,(2m)!\,(m+N-1)!\,(m+N)!\,(2N-2l+2)!}%
{4(l+m-1)(m+N-l+1)(l-1)!\,l!\,(m-1)!}\\
\times
\frac{1}{m!\,(N-l)!\,(N-l+1)!\,(2m+2N-1)!}
\HypF{4}{3}{1-l,1,1,\frac{3}{2}-l+N}%
{\frac{3}{2},2-l-m,2-l+m+N}{1}.
\end{multline}
Now we substitute $m\sim a N$ and $l\sim b N$ and perform the limit
$N\rightarrow\infty$. With Stirling's formula we determine the
limit for the quotient of factorials in front of the $_4F_3$-series
in \eqref{eq:fact4F3} as
${2\sqrt{a(a+1)}\sqrt{b(1-b)}}/({\pi(a-b+1)(a+b)})$.
For the $\HypFsimple{4}{3}$-series itself, 
we may exchange limit and summation
by uniform convergence:
\begin{equation}\notag
\lim_{N\rightarrow\infty}\HypF{4}{3}{1-l,1,1,\frac{3}{2}-l+N}%
{\frac{3}{2},2-l-m,2-l+m+N}{1}=
\HypF{2}{1}{1,1}{\frac{3}{2}}{ {(1-b)b\over(a-b+1)(a+b)} }.
\end{equation}
A combination of these results and use of the identity
(see \cite[p.~463, (133)]{Prudnikov})
\begin{equation}\notag
\HypF{2}{1}{1,1}{\frac{3}{2}}{z}=\frac{\arcsin\sqrt{z}}{\sqrt{z(1-z)}}
\end{equation}
finish the proof.
\hfill\qed

\section{Sketch of proof of Lemma~\ref{lem:complexpart}}
The method that we use for this proof is also applied 
successfully in \cite{KratBD,KratBG,KratBH,KratBI,KrZeAA}
(see in particular the tutorial description in \cite[Sec.~2]{KratBI}).

First of all, we take appropriate factors out of the determinant
in \eqref{eq:CNml2}.
To be precise, we take 
$$
\frac{(N+m-i)!}{(m+i-1)!\,(2N-2i+1)!}
$$
out of the $i$-th row of the determinant, $i=1,2,\dots,N$. Thus we obtain
\begin{multline} \label{eq:detdef}
 \prod_{i=1}^{N}\frac{(N+m-i)!}{(m+i-1)!\,(2N-2i+1)!}\\\times
	\det_{1\leq i,j\leq N}\left(
		\begin{cases}
			\po{m+i-j+1}{j-1}\po{N+j-2i+2}{N-j}\frac{N+2m-j+1}{2}
				&\text{ if } i\neq l \\
   		\po{m+i-j+1}{j-1}\po{N+j-2i+1}{N-j+1}
				&\text{ if } i=l
		\end{cases}
	\right).
\end{multline}
Let us denote the determinant in \eqref{eq:detdef} by $D(m;N,l)$.
Comparison of \eqref{eq:CNml2} and 
\eqref{eq:detdef} yields that \eqref{eq:CNml2} will be proved once
we are able to establish the determinant evaluation
\begin{multline}
D(m;N,l)=
\prod _{i=1}^{\floor{N/2}}
	\left(
		\po{m+i}{N-2i+1}\po{m+i+\frac{1}{2}}{N-2i}
	\right)
\\\times
2^{\frac{(N-1)(N-2)}{2}}\frac{
	\po{m}{N+1}\prod_{j=1}^{N}(2j-1)!
}{
	N!\prod_{i=1}^{\floor{\frac{N}{2}}}\po{2i}{2n-4i+1}
} 
\sum_{e=0}^{l-1}
(-1)^{e}\binom{N}{e}\frac{
	(N-2e)\po{\frac{1}{2}}{e}
}{
	(m+e)(m+N-e)\po{\frac{1}{2}-N}{e}
}.\label{eq:main}
\end{multline}

For the proof of \eqref{eq:main} we proceed in several steps. An outline is as
follows. In the first step we show that 
$\prod _{i=1}^{\floor{N/2}}\po{m+i}{N-2i+1}$ is a factor of $D(m;N,l)$  
as a polynomial in $m$. In the second step we
show that 
$\prod_{i=1}^{\floor{N/2}}\po{m+i+\frac{1}{2}}{N-2i}$
is a factor of $D(m;N,l)$.
In the third step we determine the maximal degree of $D(m;N,l)$ as a
polynomial in $m$, which turns out to be $\binom {N+1}2-1$. From a
combination of these three steps we are forced to conclude that
\begin{equation} \label{eq:polydef}
D(m;N,l)=\prod _{i=1}^{\floor{N/2}}
	\left(
		\po{m+i}{N-2i+1}\po{m+i+\frac{1}{2}}{N-2i}
	\right)
P(m;N,l),
\end{equation}
where $P(m;N,l)$ is a polynomial in $m$ of degree at most $N-1$.
Finally, in the fourth step, we evaluate $P(m;N,l)$ at $m=0,-1,\dots,-N$.
Namely, for $m=0,-1,\dots,-\floor{N/2}$ we show that
\begin{multline} \label{eq:polyeval}
P(m;N,l)=
(-1)^{m N+(m^2-m)/2}2^{(m^2+m)/2-N+1}\po{m}{m}
\\ \times
\frac{
		\prod_{j=1}^{N-m}(2j-1)!
\prod_{k=1}^{m}(k-1)!^2(N+k-2m-1)!\po{\frac{m-k+1}{2}}{k-1}\po{k-N}{N-m}
}{
		\prod_{i=1}^{m}(N-m-i)!(m-i)!
		\prod_{i=m+1}^{\floor{N/2}}\po{i-m}{N-2i+1}
		\prod_{i=1}^{\floor{N/2}}\po{i-m+\frac{1}{2}}{N-2i}
}.
\end{multline}
Moreover, we show that $P(m;N,l)=P(-N-m;N,l)$, which in combination with 
\eqref{eq:polyeval} gives the evaluation of $P(m;N,l)$ at $m=-\floor{N/2}-1,
\dots,-N+1,-N$.
Clearly, this determines a polynomial of maximal
degree $N-1$ uniquely. In fact, an explicit expression for $P(m;N,l)$ can 
immediately be written down using Lagrange interpolation. 
As it turns out, the resulting expression for $P(m;N,l)$ is exactly
the second line of \eqref{eq:main}.
In view of \eqref{eq:polydef},
this would establish \eqref{eq:main} and, hence, 
finish the proof of the Lemma.

\smallskip
Before going into details of these steps, it is useful to observe
the symmetry
\begin{equation}
\label{eq:simple-symmetry}
D(m;N,l)=D(m;N,N+1-l).
\end{equation}
This symmetry follows immediately from the combinatorial ``origin" of 
the determinant. For, trivially, the number of rhombus tilings which
contain the $l$-th rhombus on the symmetry axis is the same as the number
of rhombus tilings which
contain the $(N+1-l)$-th rhombus. The manipulations that finally lead to
the determinant $D(m;N,l)$ do not affect this symmetry, therefore 
$D(m;N,l)$ inherits the symmetry.

This symmetry is very useful for our considerations, because for any
claim that we want to prove (and which also obeys this symmetry) we may
freely assume $1\leq l\leq \floor{\frac{N+1}{2}}$ or
$\floor{\frac{N+1}{2}}\leq l\leq N$, whatever is more convenient.

Another useful symmetry is
\begin{equation}\label{eq:symmetry}
D(-N-m;N,l)=(-1)^{\binom {N+1}2-1}D(m;N,l).
\end{equation}
In order to establish \eqref{eq:symmetry}, we multiply the matrix
underlying $D(m;N,l)$ (as defined in \eqref{eq:detdef}) by the 
upper triangular
matrix $\left((-1)^j\binom{j-1}{i-1}\right)_{1\le i,j\le N}$. 
Using either the Gosper-Zeilberger algorithm or elementary 
``hypergeometrics" (a contiguous relation and Vandermonde summation),
the result of this multiplication is the original matrix with $m$ replaced
by $-N-m$, except that
all the entries in row $l$ have opposite sign. Hence, the equation
\eqref{eq:symmetry} follows immediately.

Now we are ready for giving details of Steps~1--4.

\bigskip
{\it Step 1. $\prod _{i=1}^{\floor{N/2}}\po{m+i}{N-2i+1}$ is a 
factor of $D(m;N,l)$}. Here, for the first time, we make use of the
symmetry \eqref{eq:simple-symmetry}. It implies, that we may restrict 
ourselves to $1\leq l\leq \floor{\frac{N+1}{2}}$.

For $i$ between $1$ and $\floor{N/2}$ 
let us consider row $N-i+1$ of the determinant $D(m;N,l)$.
Recalling that $D(m;N,l)$ is defined as the determinant in 
\eqref{eq:detdef}, we see that 
the $j$-th entry in this row has the form
$$(m+N-i-j+2)_{j-1}\,(-N+2i+j)_{N-j}\frac {N+2m-j+1} {2}.$$
Since $(-N+2i+j)_{N-j}=0$ for $j=1,2,\dots,N-2i$, the first
$N-2i$ entries in this row vanish. Therefore $(m+i)_{N-2i+1}$ is a
factor of each entry in row $N-i+1$, $i=1,2,\dots,\floor{N/2}$. Hence, the
complete product $\prod _{i=1}^{\floor{N/2}}\po{m+i}{N-2i+1}$ divides $D(m;N,l)$.

\smallskip
{\it Step 2. $\prod_{i=1}^{\floor{N/2}}\po{m+i+\frac{1}{2}}{N-2i}$
is a factor of $D(m;N,l)$}. Again we make use of the symmetry 
\eqref{eq:simple-symmetry}, which allows us to restrict 
ourselves to $1\leq l\leq \floor{\frac{N+1}{2}}$. 

We observe that the product can be rewritten as
$$\prod_{i=1}^{\floor{N/2}}\po{m+i+\frac{1}{2}}{N-2i}=
\prod _{e=1} ^{N-2}(m+e+\frac {1} {2})^{\min\{e,N-e-1\}}.$$
Therefore, 
because of the other symmetry \eqref{eq:symmetry}, it suffices to
prove that $(m+e+1/2)^{e}$ divides $D(m;N,l)$ for $e=1,2,\dots,\floor{N/2}-1$.
In order to do so, we claim that for each such $e$ there are $e$ 
linear combinations of the columns, which are themselves linearly
independent, that vanish for $m=-e-1/2$. More precisely, we claim that
for $k=1,2,\dots,e$ there holds
\begin{multline} \label{eq:lincomb}
\sum _{j=1} ^{k}\binom kj \cdot(\text {column
$(N+1-2e+k+j)$ of $D(-e-1/2;N,l)$})\\
-\frac{\po{N-e-l+1}{k}}{(-4)^{k}\po{N-e-l+\frac{3}{2}}{k}}
\cdot(\text {column
$(N+1-2e)$ of $D(-e-1/2;N,l)$})=0.
\end{multline}
As is not very difficult to see (cf\@. \cite[Sec.~2]{KratBI}) this would
imply that $(m+e+1/2)^e$ divides $D(m;N,l)$.

Obviously, a proof of \eqref{eq:lincomb} amounts to proving two
hypergeometric identities, one for the restriction of \eqref{eq:lincomb}
to the $i$-th row, $i\ne l$, and another for 
the restriction of \eqref{eq:lincomb}
to the $l$-th row itself. Both identities can be easily established by
using either the Gosper--Zeilberger algorithm or elementary
``hypergeometrics" (again, a contiguous relation and Vandermonde summation).

\smallskip
{\it Step 3. $D(m;N,l)$ is a polynomial in $m$ of maximal degree
$\binom {N+1}2-1$}.
Clearly, the degree in $m$ of the $(i,j)$-entry in
the determinant $D(m;N,l)$ is $j$ for $i\ne l$, while it is $j-1$ for
$i=l$. Hence, in the defining expansion of the determinant, each term
has degree $\left(\sum _{j=1} ^{N}j\right)-1=\binom {N+1}2-1$.

\smallskip
{\it Step 4. Evaluation of $P(m;N,l)$ at $m=0,-1,\dots,-N$}.
This step is the most technical one, therefore we shall be only brief here.

Again, we make use of the symmetry \eqref{eq:simple-symmetry}, 
and this time restrict ourselves to $\floor{\frac{N+1}{2}}\leq l\leq N$.
On the other hand, by the symmetry \eqref{eq:symmetry} and by
the definition \eqref{eq:polydef} of $P(m;N,l)$, we have
$P(m;N,l)=P(-N-m;N,l)$. Therefore, it suffices to
compute the evaluation of $P(m;N,l)$ at $m=0,-1,\dots,-\floor{N/2}$.

What we would like to do is, for any $e$ with $0\le e\le \floor{N/2}$,
to set $m=-e$ in
\eqref{eq:polydef}, compute $D(-e;N,l)$, and then express $P(-e;N,l)$ as
the ratio of $D(-e;N,l)$ and the right-hand side product evaluated
at $m=-e$. Unfortunately, this is typically a ratio $0/0$ and, hence,
undetermined. So, we have to first divide both sides of \eqref{eq:polydef}
by the appropriate power of $(m+e)$, and only then set $m=-e$.

Let $e$, $0\le e\le \floor{N/2}$, be fixed.
For $k=0,1,\dots,{e-1}$ we add 
\begin{equation}
\label{eq:colOp}
\sum_{i=1}^k\binom{k}{i}\cdot(\text {column
$(N+1-2e+k+i)$ of $D(m;N,l)$})
\end{equation}
to column ${N+1-2e+k}$ of $D(m;N,l)$. 
The effect (which is again proved by
either the Gosper--Zeilberger algorithm or ``hypergeometrics") 
is that then $(m+e)$ is a factor
of each entry in column ${N+1-2e+k}$. So, we take $(m+e)$ out of
each entry of column ${N+1-2e+k}$, $k=0,1,\dots,{e-1}$. 

Let
$D_2(m;N,l,e)$ denote the resulting determinant. From what we did so far,
it is straight-forward that we must have
$$D(m;N,l)=(m+e)^e D_2(m;N,l,e).$$
A combination with \eqref{eq:polydef} gives that
\begin{equation} \label{eq:polydet}
P(m;N,l)=D_2(m;N,l,e)\prod _{i=1}^{\floor{N/2}}
	\left(
		\po{m+i}{e-i}\po{m+e+1}{N-i-e}\po{m+i+\frac{1}{2}}{N-2i}
	\right)^{-1}.
\end{equation}
Now, in order to determine the evaluation of $P(m;N,l)$ at $m=-e$,
we set $m=-e$ in \eqref{eq:polydet}. It turns out that 
the determinant $D_2(-e;N,l,e)$ vanishes for $e\ge N+1-l$, whereas
for $e<N+1-l$ the matrix underlying
$D_2(-e;N,l,e)$ has a block form as illustrated in Figure~\ref{fig:matrix}. 
Therefore, in the latter case, the determinant $D_2(-e;N,l,e)$ equals the product of the
determinants of $Q_1$, $Q_2$, and $M$, each of which can be easily
evaluated explicitly. For, $Q_1$ and $Q_2$ are upper and lower triangular
matrices,
respectively, and the determinant of $M$ is again easily determined by
applying the determinant lemma \cite[Lemma~2.2]{KratAM}. Then, by
combining these computations with \eqref{eq:polydet},
and performing some simplification, the evaluation \eqref{eq:polyeval}
follows.

\begin{figure}
	\begin{picture}(150,190)(-150,-20)
		\put(0,0){\line(0,1){150}}
		\put(0,150){\line(1,0){150}}
		\put(150,150){\line(0,-1){150}}
		\put(150,0){\line(-1,0){150}}
		\put(50,0){\line(0,1){100}}
		\put(50,145){\line(0,1){5}}
		\put(100,0){\line(0,1){150}}
		\put(0,50){\line(1,0){100}}
		\put(145,50){\line(1,0){5}}
		\put(0,100){\line(1,0){150}}
		\put(0,100){\line(1,-1){100}}
		\put(122,45){{\Large 0}}
		\put(22,20){{\Large 0}}
		\put(30,80){{\Large 0}}
		\put(65,15){{\Large 0}}
		\put(71,71){{\Huge $*$}}
		\put(46,121){{\Huge $*$}}
		\put(120,120){{\Large $M$}}
		\put(13,63){{\Large $Q_2$}}
		\put(79,31){{\Large $Q_1$}}
		\put(2,-7){\tiny $1$}
		\put(52,-7){\tiny $N+1-2e$}
		\put(102,-7){\tiny $N+1-e$}
		\put(28,155){\tiny $N-2e$}
		\put(83,155){\tiny $N-e$}
		\put(145,155){\tiny $N$}
		\put(-5,145){\tiny $1$}
		\put(-19,95){\tiny $e+1$}
		\put(-36,45){\tiny $N+1-e$}
		\put(155,103){\tiny $e$}
		\put(155,53){\tiny $N-e$}
		\put(155,3){\tiny $N$}
	\end{picture}
\caption{}
\label{fig:matrix}
\end{figure}

\smallskip
This finishes the proof of Lemma~\ref{lem:complexpart}.
\hfill\qed

\ifx\undefined\bysame
\newcommand{\bysame}{\leavevmode\hbox to3em{\hrulefill}\,}
\fi

\end{document}